\documentclass[11pt]{amsart}

\usepackage{amssymb}
\usepackage[frenchb]{babel}

\evensidemargin\oddsidemargin

\begin{document}
\baselineskip=18pt
\setcounter{page}{1}
    
\newtheorem{Theo}{Th\'eor\`eme\!\!}
\newtheorem{Rqs}{Remarque\!\!}
\newtheorem{Lemm}{Lemme}
\newtheorem{Corru}{Corollaire 1\!\!}
\newtheorem{Corrd}{Corollaire 2\!\!}
\newtheorem{Rq}{Remarque\!\!}

\renewcommand{\theRqs}{}
\renewcommand{\theTheo}{}
\renewcommand{\theCorrd}{}
\renewcommand{\theCorru}{}
\renewcommand{\theRq}{}

\def\a{\alpha}
\def\b{\beta}
\def\B{{\bf B}} 
\def\C{{\mathcal{C}}} 
\def\CC{{\mathbb{C}}} 
\def\E{{\mathcal{E}}} 
\def\Ea{E_\a}
\def\EE{{\mathbb{E}}} 
\def\elaw{\stackrel{d}{=}}
\def\eps{\varepsilon}
\def\F{{\bf F}} 
\def\G{\gamma} 
\def\G{{\bf \Gamma}} 
\def\HH{{\mathbb{H}}} 
\def\hS{{\hat S}}
\def\hT{{\hat T}}
\def\hX{{\hat X}}
\def\ii{{\rm i}}
\def\K{{\bf K}} 
\def\L{{\mathcal{L}}} 
\def\lb{\lambda}
\def\lacc{\left\{}
\def\lcr{\left[}
\def\lpa{\left(}
\def\lva{\left|}
\def\NN{{\mathbb{N}}} 
\def\pb{{\mathbb{P}}}
\def\R{{\mathcal{R}}}
\def\rl{{\mathbb{R}}}
\def\racc{\right\}}
\def\rcr{\right]}
\def\rpa{\right)}
\def\rva{\right|}
\def\T{{\bf T}} 
\def\TT{{\rm T}} 
\def\U{{\bf U}} 
\def\Un{{\bf 1}}
\def\ZZ{{\mathbb{Z}}} 
      
\def\abstractname{R\'esum\'e}

\renewcommand{\refname}{Bibliographie}

\def\keywordsname{Mots-cl\'es}

\def\subjclassname{Classification AMS}

\bibliographystyle{french}

\newcommand{\fin}{\vspace{-0.4cm}
                  \begin{flushright}
                  \mbox{$\Box$}
                  \end{flushright}
                  \noindent}

\title[Produit Beta-Gamma et r\'egularit\'e du signe]{Produit Beta-Gamma et r\'egularit\'e du signe}

\author[Thomas Simon]{Thomas Simon}

\address{Laboratoire Paul Painlev\'e, Universit\'e Lille 1, F-59655 Villeneuve d'Ascq Cedex. {\em Adresse \'electronique} : {\tt simon@math.univ-lille1.fr}}



\begin{abstract} 
On \'etudie la positivit\'e compl\`ete du noyau de convolution multiplicatif T associ\'e au produit de deux variables al\'eatoires ind\'ependantes $\B(a,b)$ et $\G(c).$ Ce noyau T est compl\`etement positif d'ordre infini si $b\in\NN^*$ ou si $d = a+b -c\in\NN.$ Dans les autres cas la r\'egularit\'e du signe de T a toujours un ordre fini, qui est ici calcul\'e. Plus pr\'ecis\'ement, pour tout $n\ge 1$ on montre que T est compl\`etement positif d'ordre $n + 1$ si et seulement si $(d,b)$ est situ\'e au dessus d'un certain escalier $\E_n$ dessin\'e dans le demi-plan sup\'erieur. Cet escalier caract\'erise aussi la constance du signe de plusieurs d\'eterminants associ\'es \`a la fonction hyperg\'eom\'etrique confluente de seconde esp\`ece.
\end{abstract}

\maketitle

\section{Introduction}

Soit $I$ un intervalle de $\rl$ et $K$ un noyau \`a valeurs r\'eelles d\'efini sur $I\times I.$ On dit que $K$ est compl\`etement positif d'ordre $n$ (${\rm TP}_n$ selon l'acronyme anglo-am\'ericain) si 
$$\det \lcr K(x_i, y_j)\rcr_{1\le i,j\le m}\; \ge \; 0$$ 
pour tout $m\in [1,n], x_1< \ldots < x_m$ et $y_1< \ldots < y_m.$  Quand la propri\'et\'e est vraie pour tout $n$ on dit que $K$ est compl\`etement positif d'ordre infini (${\rm TP}_\infty$). On dit que $K$ est r\'egulier pour le signe \`a l'ordre $n$ (${\rm SR}_n$) s'il existe $\{\eps_i\}_{1\le i\le n}$ dans $\{-1,1\}$ tel que
$$\eps_m \det \lcr K(x_i, y_j)\rcr_{1\le i,j\le m}\; \ge \; 0$$ 
pour tout $m\in [1,n], x_1< \ldots < x_m$ et $y_1< \ldots < y_m.$ Quand la propri\'et\'e est vraie pour tout $n$ on dit que $K$ est r\'egulier pour le signe \`a l'ordre infini (${\rm SR}_\infty$). Ces quatre propri\'et\'es sont dites strictes lorsque les in\'egalit\'es correspondantes sont toutes strictes, et on utilise alors la notation ${\rm STP}_n, {\rm SSR}_n, {\rm STP}_\infty$ ou ${\rm SSR}_\infty.$ Introduites avant-guerre par Krein, Schoenberg et Gantmacher, ces notions ont ensuite connu un d\'eveloppement important qui a abouti \`a l'ouvrage classique \cite{K}. On renvoie \`a la monographie \cite{P} pour des avanc\'ees plus r\'ecentes ainsi qu'un historique du sujet.

Une fonction $f : \rl \to\rl^+$ est dite fr\'equence de P\'olya d'ordre $n\le \infty$ (${\rm PF}_n$) si le noyau $K(x,y) = f(x-y)$ est ${\rm TP}_n$ sur $\rl\times\rl.$ Les densit\'es de probabilit\'e dans la classe ${\rm PF}_\infty$ ont \'et\'e caract\'eris\'ees par Schoenberg - voir le th\'eor\`eme 7.3.2 (a) p. 345 dans \cite{K} - par leur transform\'ee de Laplace, dont l'inverse se factorise sous la forme
\begin{equation}
\label{droite}
\frac{1}{\EE[e^{sX}]} \; =\; e^{-\gamma s^2+\delta s}\prod_{n=0}^\infty (1+a_n s)e^{-a_n s}
\end{equation}
o\`u $X$ est la variable al\'eatoire associ\'ee, avec $\gamma \ge 0, \delta\in\rl,$ et $\sum a_n^2 < \infty.$ L'exemple classique est la densit\'e gaussienne. Le cas des densit\'es sur la demi-droite a \'egalement \'et\'e caract\'eris\'e par une factorisation analogue  - voir le th\'eor\`eme 7.3.2 (b) dans \cite{K} - qui s'\'ecrit 
\begin{equation}
\label{ddroite}
\frac{1}{\EE[e^{-sX}]} \; =\; e^{\delta s}\prod_{n=0}^\infty (1+a_n s)
\end{equation}
avec $\delta, a_n \ge 0$ et $\sum a_n < \infty.$ Ceci montre que les densit\'es ${\rm PF}_\infty$ sur $\rl^+$ sont toutes, modulo une d\'erive, des convol\'ees de densit\'es exponentielles de param\`etres diff\'erents et forment donc une sous-classe des densit\'es dites GGC - voir \cite{Bo}. Les fonctions ${\rm PF}_2$ sont ais\'ement caract\'eris\'ees par la log-concavit\'e sur leur support - voir le th\'eor\`eme 4.1.9 dans \cite{K}, ce qui a pour cons\'equence utile la stabilit\'e de la propri\'et\'e ${\rm PF}_2$ par multiplication ponctuelle. Une telle propri\'et\'e n'est plus vraie pour les classes ${\rm PF}_n, n\ge 3,$ qui sont en g\'en\'eral plus difficiles \`a \'etudier que les classes ${\rm PF}_2$ ou ${\rm PF}_\infty.$ Toutes les propri\'et\'es ${\rm PF}_n$ sont cependant stables par produit de convolution - voir la proposition 7.1.5 dans \cite{K}, ce qui dans le cas $n=2$ est un cas particulier du th\'eor\`eme de Pr\'ekopa, et dans le cas des densit\'es signifie que la somme ind\'ependante de deux variables al\'eatoires dont les densit\'es sont ${\rm PF}_n$ a encore une densit\'e ${\rm PF}_n.$ Cette stabilit\'e par produit de convolution sera importante dans la suite de cet article.

La notion de fonction fr\'equence de P\'olya est donc associ\'ee de fa\c{c}on naturelle \`a l'addition des variables al\'eatoires. Pourtant, les exemples les plus standard de noyaux compl\`etement positifs sont li\'es \`a la convolution multiplicative. Ainsi, le caract\`ere ${\rm TP}_\infty$ du noyau $e^{xy}$ mentionn\'e au paragraphe 1.2.(i) p. 15 dans \cite{K} est \'equivalent \`a la m\^eme propri\'et\'e pour le noyau $f (xy^{-1}),$ o\`u
$$f(x) \; =\; \frac{x^{c-1} e^{-x}}{\Gamma(c)} \Un_{(0,\infty)}(x)$$
est la densit\'e de la variable $\G(c).$ Cette propri\'et\'e qui est \'equivalente au caract\`ere ${\rm PF}_\infty$ de $\log\G(c),$ peut aussi  se voir par (\ref{droite}) en \'ecrivant les moments fractionnaires de $\G(c)$ :
$$\frac{1}{\EE[\G(c)^s]}\; =\; \frac{\Gamma(c)}{\Gamma(c+s)}\; =\; e^{\gamma^{}_c s}\prod_{n=0}^{\infty} \lpa 1 +s(n+c)^{-1}\rpa e^{-s(n+c)^{-1}}$$
pour un certain $\gamma^{}_c > 0,$ o\`u la deuxi\`eme \'egalit\'e vient de la formule de Weierstrass - voir par exemple \cite{E} p. 2. De m\^eme le caract\`ere ${\rm TP}_\infty$ du noyau $\Un_{\{0\le x\le y\le 1\}}$ mentionn\'e au paragraphe 1.2.(ii) de \cite{K}, est pr\'ecis\'ement celui du noyau de convolution multiplicatif de la variable $\U$ uniforme sur $(0,1).$ Ceci signifie que $-\log \U\sim$ Exp($1$) a une densit\'e ${\rm PF}_\infty,$ et est bien s\^ur une cons\'equence \'evidente de (\ref{ddroite}) et de l'expression des moments fractionnaires de $\U$.

Le livre \cite{K} fournit plusieurs autres exemples de noyaux compl\`etement positifs reli\'es \`a des variables al\'eatoires, souvent dans le cadre des familles exponentielles. Le cas du noyau multiplicatif associ\'e \`a la variable $\B(a,b)$ de densit\'e
$$\frac{\Gamma(a+b)}{\Gamma(a)\Gamma(b)} \, x^{a-1} (1-x)^{b-1} \Un_{(0,1)}(x)$$
est aussi implicitement trait\'e par le th\'eor\`eme 10.1.1 p. 503, pour $b\in\NN^*$. En effet, en choisissant $w_1\equiv\ldots\equiv w_n\equiv 1$ dans (1.4) p. 502, ce th\'eor\`eme montre le caract\`ere ${\rm TP}_\infty$ du noyau $(x-y)_+^n$ et donc du noyau multiplicatif associ\'e \`a $\B(a,n+1).$ Ce dernier fait, qui remonte \`a Schoenberg et Whitney en 1953, se voit aussi \`a l'aide de l'identit\'e
\begin{equation}
\label{BetaU}\B(a,n+1) \; \elaw \; \U^{1/a}\;\times\;\cdots\;\times\; \U^{1/(a+n)}
\end{equation}
obtenue en comparant les moments fractionnaires. Le cas plus d\'elicat $b\notin\NN^*$ est une cons\'equence du th\'eor\`eme 3.4(i) dans \cite{DM}, o\`u ce probl\`eme semble avoir \'et\'e consid\'er\'e pour la premi\`ere fois.

\begin{Theo}[Dette-Munk] Le noyau de convolution multiplicatif associ\'e \`a $\B(a,b)$ est ${\rm TP}_{[b+1]}$ pour tous $a,b >0.$
\end{Theo} 

Pour obtenir l'\'enonc\'e ci-dessus, on a fait la correction $n\to n+1$ dans le r\'esultat de Dette et Munk, ce qui se voit ais\'ement en appliquant correctement le  lemme 3.1 dans \cite{DM}. Ce m\^eme lemme, et des arguments classiques de \cite{K} que nous ne d\'etaillerons pas ici - voir la remarque (c) \`a la fin de la troisi\`eme partie du pr\'esent article, montrent que ce m\^eme noyau n'est en revanche pas ${\rm TP}_{[b+2]}$ ni d'ailleurs ${\rm SR}_{[b+2]}$, sauf pour $b\in\NN^*.$ Quand $b\notin\NN^*,$ la r\'egularit\'e du signe du noyau multiplicatif de $\B(a,b)$ est ainsi une fonction croissante de $[b],$ et ceci peut se deviner par l'identit\'e  
\begin{equation}
\label{crs}
\B(a,b+1) \; \elaw \; \B(a,b)\;\times\; \U^{1/(a+b)}
\end{equation}
puisque $\U$ - et donc $\U^{1/(a+b)}$ de mani\`ere \'equivalente - est associ\'ee \`a un noyau ${\rm TP}_\infty.$ D'une mani\`ere g\'en\'erale, on s'attend \`a ce que le produit ind\'ependant par des variables ayant une telle propri\'et\'e ${\rm TP}_\infty$ r\'egularise le signe. 

Le but de cet article est d'\'etudier ce ph\'enom\`ene pour les variables produit 
$$\T(a,b,c)\; =\; \B(a,b)\times \G(c),$$ 
dont la densit\'e sur $(0,+\infty)$ s'\'ecrit
$$\frac{\Gamma (a +b) x^{c-1}e^{-x}}{\Gamma
  (a)\Gamma(b)\Gamma(c)}\int_0^\infty e^{-tx} \frac{t^{b-1}}{(t +1)^{d}} dt\; =\; \frac{\Gamma (a +b) x^{c-1}}{\Gamma
  (a)\Gamma(c)}\; e^{-x}\Psi (b, b+1-d, x)$$
avec la notation $d = a+b - c.$ Ci-dessus, $\Psi$ d\'esigne la fonction hyperg\'eom\'etrique confluente de seconde esp\`ece, encore appel\'ee fonction de Tricomi, pour les d\'etails de laquelle on renvoie au chapitre 6 de \cite{E}. Par le th\'eor\`eme 1.2.1 dans \cite{K}, la r\'egularit\'e du signe du noyau de convolution multiplicatif associ\'e \`a $\T(a,b,c)$ est la m\^eme que celle du noyau
$$\TT_{b,d} (x,y) \; =\; e^{-xy^{-1}}\Psi (b, b+1-d, xy^{-1})$$
sur $(0,+\infty)\times(0,+\infty),$ et ne d\'epend donc que du couple $(b,d).$ On sait par le r\'esultat pr\'ec\'edent que cette r\'egularit\'e est croissante en $[b]$ et on s'attend \`a ce qu'elle le soit aussi en $[d]$ en vertu de l'identit\'e
\begin{equation}
\label{crescent}
\T(a,b,c)\; \elaw\;\U^{1/c}\;\times\; \T(a,b,c+1).
\end{equation}
Introduisons maintenant nos notations. Sur le demi-plan sup\'erieur $\HH = \{y\ge 0\}$ on consid\`ere le r\'eseau $\R = \{x\in\NN\}\cup\{y\in\NN^*\}$ et, pour tout $n\ge 1,$ l'escalier $\E_n$ constitu\'e des $(n-1)$ marches et  $n$ contremarches 
$$\bigcup_{k=1}^{n-1}\{n-k-1\le x\le n-k, y=k\}\quad\mbox{et}\quad\bigcup_{k=0}^{n-1}\{x=k, n-k-1\le y\le n-k\},$$
et des deux marches infinies $\{x\le 0, y=n\}$ et $\{x\ge n-1, y=0\}.$ On consid\`ere sur $\HH$ la relation $(x,y)\succ(z,t)\Leftrightarrow x\ge z$ ou $y\ge t.$ Pour tout $A\subset\HH$ et $(x,y)\in\HH$ on note $(x,y)\succ A$ si $\forall\,(z,t)\in A$ on a $(x,y)\succ (z,t).$ Par abus d'\'ecriture on note $(x,y)\prec A$ si $\exists\,(z,t)\in A$ tel que $x < z$ et $y<t.$ Visuellement parlant, pour tout $n\ge 1$ on a $(x,y)\prec\E_n$ quand $(x,y)$ est sous l'escalier et $(x,y)\succ\E_n$ quand  $(x,y)$ est au dessus. Notre r\'esultat principal s'\'enonce comme suit.

\begin{Theo} {\rm (i)} On a $\TT_{b,d}\in{\rm TP}_\infty\Leftrightarrow (d,b)\in\R.$ 

\noindent
{\rm (ii)} Si $(d,b)\notin\R,$ alors $\forall\,n\ge 1$ on a $(d,b)\succ\E_n\,\Leftrightarrow\,\TT_{b,d}\in{\rm TP}_{n+1}$ et $(d,b)\prec\E_n\,\Rightarrow\,\TT_{b,d}\notin{\rm SR}_{n+1}.$
\end{Theo}

Ce r\'esultat montre en particulier que la variable $\G(c)$ r\'egularise le signe dans le produit $\G(c)\times\B(a,b)$ seulement pour $c\le a+b.$ Dans le quadrant positif $\{b > 0, d > 0\},$ il met aussi en \'evidence une grille de croissance de la propri\'et\'e ${\rm TP}_n$ pour $\T(a,b,c),$ dont l'allure n'est pas sans rappeler celle du nombre de racines positives de $\Psi$ pour les valeurs n\'egatives de ses param\`etres - voir la figure p. 290 dans \cite{E}. Malgr\'e le r\^ole central de $\Psi$ dans notre argument et les diverses manipulations dont cette fonction peut \^etre l'objet, l'analogie entre les deux grilles est sans doute fortuite.

La caract\'erisation de la propri\'et\'e ${\rm TP}_\infty,$ qui d\'ecoule sans difficult\'e du th\'eor\`eme de Schoenberg et d'identit\'es en loi simples, est  d\'emontr\'ee dans la deuxi\`eme partie de cet article o\`u l'on donne aussi une nouvelle  preuve de la transformation de Kummer pour $\Psi$, laquelle fait appara\^\i tre une sym\'etrie dans le quadrant positif pour la caract\'erisation de la propri\'et\'e ${\rm TP}_{n+1}$. Dans la troisi\`eme partie on d\'emontre le deuxi\`eme volet du th\'eor\`eme, qui est bien moins facile que le premier. On ram\`ene d'abord le probl\`eme \`a la constance du signe d'un certain d\'eterminant de Hankel avec poids. On montre ensuite le r\'esultat par divers arguments m\^elant alg\`ebre lin\'eaire, analyse asymptotique et l'identit\'e (\ref{crescent}). Cette derni\`ere est cruciale pour r\'eduire le domaine d'\'etude \`a deux s\'eries de marches adjacentes situ\'ees de part et d'autre de $\E_n.$ 

Des d\'eterminants fonctionnels identiques aux n\^otres ont \'et\'e r\'ecemment \'etudi\'es dans \cite{BI, IL} afin d'\'etablir de nouvelles in\'egalit\'es de type Tur\'an. Notre r\'esultat principal permet de pr\'eciser  certaines de ces in\'egalit\'es, ce que nous faisons dans la quatri\`eme partie. Comme dans l'article classique \cite{KS} nous nous int\'eressons aussi \`a des Wronskiens pour des familles de fonctions $\Psi$ ayant des param\`etres contigus, dont nous caract\'erisons la constance du signe par les fronti\`eres $\E_n.$ La forme en escalier de ces derni\`eres est peut-\^etre un peu surprenante, et je n'ai pas trouv\'e d'exemples de fronti\`ere de ce type dans la litt\'erature ancienne ou moderne consacr\'ee aux d\'eterminants fonctionnels.

\section{Quelques r\'esultats pr\'eliminaires}

\subsection{Le cas ${\rm TP}_\infty$} Dans ce paragraphe on consid\`ere la partie (i) du th\'eor\`eme, qui se montre facilement gr\^ace \`a la caract\'erisation de Schoenberg. Les moments fractionnaires de $\T(a,b,c)$ sont donn\'es par 
\begin{equation}
\label{FracT}
\EE[\T(a,b,c)^s]\; =\; \frac{\Gamma(a+s)\Gamma(c+s)\Gamma(a+b)}{\Gamma(a)\Gamma(c)\Gamma(a+b+s)}
\end{equation}
pour tout $s\in\CC$ tel que $\Re(s) > (-a)\vee(-c).$ On en d\'eduit les identit\'es 
$$\T(a,b, a+b)\;\elaw \; \G(a)\qquad\mbox{et}\qquad
\T(a,b,a+b+n) \; \elaw \;\G(a)\,\times\,\U^{1/a}\,\times\,\cdots\,\times\, \U^{1/(a+n-1)}$$
pour tout entier $n\ge 1$ (la premi\`ere est un des nombreux exemples de formules d\'efinissant l'alg\`ebre Beta-Gamma \cite{D}). Ces identit\'es montrent que $\TT_{b,d}$ est ${\rm TP}_\infty$ si $d\in\NN.$ En vertu de (\ref{BetaU}) on a d'autre part 
$$\T(a,n,c)\;\elaw\;\U^{1/a}\;\times\;\cdots\;\times\; \U^{1/(a+n-1)}\;\times\;\G(c)$$
pour tout entier $n\ge 1,$ et ceci montre que $\TT_{b,d}$ est ${\rm TP}_\infty$ si $b\in\NN^*.$ Si $(d,b)\notin\R,$ on voit enfin que la fonction
$$s\;\mapsto\;\frac{1}{\EE[\T(a,b,c)^s]}$$
admet les nombres $\{-(a+b+n), \,n\in\NN\}$ pour p\^oles simples et n'est donc pas holomorphe sur $\CC.$ A fortiori elle n'est pas du type d\'efini dans (\ref{droite}), et le th\'eor\`eme 3.2 (a) p. 345 dans \cite{K} entra\^\i ne que la densit\'e de $\log \T(a,b,c)$ n'est pas ${\rm PF}_\infty$, autrement dit $\TT_{b,d}$ n'est pas ${\rm TP}_\infty.$ On a donc bien d\'emontr\'e que
$$\TT_{b,d}\in {\rm TP}_\infty\;\Leftrightarrow\; (d,b)\in\R.$$

\subsection{La transformation de Kummer et deux applications} Supposons maintenant $c < a+b,$ autrement dit $d > 0.$ La formule (\ref{FracT}) entra\^\i ne l'identit\'e 
$$\T(a,b,c) \; \elaw \; \T(c,d,a)$$
et en comparant les densit\'es de part et d'autre on trouve
\begin{eqnarray}
\label{Kumm}
\Psi(b, b+1-d,x) \; =\; x^{d-b}\Psi (d, d+1-b, x)
\end{eqnarray}
pour tout $b,d,x > 0.$ Cette \'egalit\'e est connue sous le nom de transformation de Kummer pour la fonction de Tricomi et est obtenue dans le chapitre 6 de \cite{E} pp. 256-257 en comparant diff\'erentes repr\'esentations de Mellin-Barnes de $\Psi.$ Cette derni\`ere m\'ethode est plus elabor\'ee que l'argument pr\'ec\'edent reposant sur une b\^ete identification de moments fractionnaires, mais son avantage est de pouvoir \'etendre la transformation \`a des valeurs n\'egatives de $b$ et $d$ - voir le haut de la page 257 dans le chapitre 6 de \cite{E}. On renvoie \`a \cite{D} pour d'autres relations entre repr\'esentations de Mellin-Barnes et identit\'es en loi multiplicatives.\\

La formule (\ref{Kumm}) entra\^\i ne que
$$\TT_{b,d}\in {\rm TP}_{n+1}\;\Leftrightarrow\; \TT_{d,b}\in {\rm TP}_{n+1}$$
pour tout $b,d > 0$ et $n\ge 1.$ Cette sym\'etrie par rapport \`a la premi\`ere bissettrice sur le quadrant positif est bien s\^ur compatible avec celle de la fronti\`ere $\E_n$. Combin\'ee avec le r\'esultat de Dette et Munk, elle montre que
$$\TT_{b,d}\in {\rm TP}_{n+1}\quad\mbox{si}\quad \inf(b,d)\ge n.$$
On verra par la suite qu'on peut aussi obtenir ce r\'esultat directement en consid\'erant un certain d\'eterminant de Hankel associ\'e \`a la fonction $\Psi.$ \\

On donne maintenant, en passant, une deuxi\`eme application de la transformation de Kummer ind\'ependante de tout le reste de cet article. Soit la variable al\'eatoire positive $\K^2(a,b,c)$ ayant pour densit\'e
$$\frac{x^{a-1} (1+x)^{-(a+b)}e^{-cx}}{\Gamma(a)\Psi(a,1-b,c)} \Un_{(0,\infty)}(x)$$
o\`u $a,c > 0$ et $b\in\rl.$ Cette variable appara\^\i t en statistique bay\'esienne de r\'eseaux de file d'attente o\`u elle porte le nom de variable de Kummer de type 2. On suppose maintenant $b > 0$ et on calcule la transform\'ee de Laplace de $\K^2(a+b,-b,c)$ :
\begin{eqnarray*}
\EE\lcr e^{-\lb \K^2(a+b,-b,c)} \rcr & = & \frac{\Psi(a+b,1+b,c+\lb)}{\Psi(a+b,1+b,c)}\\
& = & \frac{\Psi(a,1-b,c+\lb)}{\Psi(a,1-b,c)}\lpa \frac{c}{c+\lb}\rpa^b
\; = \; \EE\lcr e^{-\lb \K^2(a,b,c)}\rcr\EE\lcr e^{-\lb c^{-1} \G(b)}\rcr,
\end{eqnarray*}
o\`u la deuxi\`eme \'egalit\'e vient imm\'ediatement de (\ref{Kumm}). On en d\'eduit l'\'egalit\'e en loi additive
$$\K^2(a+b,-b,c)\; =\; \K^2(a,b,c)\, +\, c^{-1} \G(b)$$
pour tous $a,b,c >0,$ r\'ecemment obtenue par Koudou et Vallois - voir (2.24) dans \cite{KV} - comme cons\'equence d'une certaine propri\'et\'e d'ind\'ependance caract\'erisant le couple $(\K^2(a,b,c), \G(b)).$ Le calcul direct des transform\'ees de Laplace, plus imm\'ediat, avait d\'ej\`a \'et\'e fait par Letac - voir le th\'eor\`eme 4.1 de \cite{L} - avec un argument un peu diff\'erent du n\^otre utilisant sur la fonction hyperg\'eom\'etrique confluente de premi\`ere esp\`ece - voir la formule 6.6(7) dans \cite{E}. 

\subsection{R\'esultats d'alg\`ebre lin\'eaire} On \'enonce dans ce paragraphe quatre lemmes \'evaluant des d\'eterminants, qui serviront dans la r\'eciproque de la partie (ii) du th\'eor\`eme. On utilise la notation $0^0 = 1$ et $\det_n[a_{ij}] = \det [a_{ij}]_{0\le i,j\le n}$ pour tout tableau de nombres $\{ a_{ij}\}$ et tout $n\in\NN.$ On rappelle aussi l'\'ecriture usuelle du symbole de Pochammer : $(z)_0 = 1$ et $(z)_k = z(z+1)\ldots(z+k-1).$ Enfin, le dernier lemme reprend les m\^emes notations matricielles que celles du chapitre 0 de \cite{K}.

\begin{Lemm} 
\label{Det1}
On a ${\det}^{}_n \lcr (d-n+j)_i\rcr\, =\, \prod_{k=0}^n k!$ pour tout $d\in\rl$ et $n\in\NN.$
\end{Lemm}

\noindent
{\em Preuve} : Chaque ligne d'indice $i$ du d\'eterminant est une fonction polynomiale en $d$ de degr\'e $i$, dont la d\'eriv\'ee est une combinaison lin\'eaire \`a coefficients constants des lignes d'indice $0,\ldots, i-1.$ Ceci entra\^\i ne que $d\mapsto {\det}^{}_n \lcr (d-n+j)_i\rcr$ est constante, et vaut donc ${\det}^{}_n \lcr (j-n)_i\rcr = \prod_{k=0}^n k!$ comme souhait\'e.\fin

\begin{Rqs}{\em (a) Par des op\'erations \'el\'ementaires sur les colonnes, on peut aussi voir directement la relation de r\'ecurrence ${\det}^{}_n \lcr (d-n+j)_i\rcr\, =\,n!\times{\det}^{}_{n-1} \lcr (d-n+1+j)_i\rcr$ et conclure avec la condition initiale ${\det}^{}_0 \lcr (d+j)_i\rcr = 1.$

(b) Ce d\'eterminant a quelque cousinage avec celui \'etudi\'e dans \cite{KS} p. 31, qui s'\'ecrit  ${\det}^{}_n \lcr \{d+j\}_i\rcr$ avec la notation $\{z\}_0 = 1$ et $\{z\}_k = z(z-1)\ldots(z+1-k),$ et qui est aussi un premier terme dans un d\'eveloppement asymptotique. Ce dernier est cependant un Vandermonde et d\'epend de $d.$ }
\end{Rqs}

\begin{Lemm} 
\label{Det2}
Pour tout entier $q\ge 0$ et tous r\'eels $\mu,\nu$ tels que $\mu > \nu+q,$ on a 
$${\det}^{}_q \lcr (\nu-q+j)_i\Gamma(\mu-\nu +q- i-j)\rcr\, =\, \prod_{k=0}^q k!\,\Gamma (\mu-\nu -k) (\mu-k-1)^{q-k}.$$
\end{Lemm}

\noindent
{\em Preuve} : On commence par mettre $\prod_{k=0}^q \Gamma (\mu-\nu -k)$ en facteur puis, en ajoutant successivement la $(q-i)$-i\`eme ligne \`a la $(q-i+1)$-i\`eme pour $i = 1, \ldots, q$ on factorise par $(\mu-1)^q.$ On r\'ep\`ete l'op\'eration pour $i=1, \ldots, q-1$ et on factorise par $(\mu-2)^{q-1}...$ Ainsi de suite jusqu'\`a trouver 
$${\det}^{}_q \lcr (\nu-q+j)_i\Gamma(\mu-\nu +q- i-j)\rcr\; =\; {\det}^{}_q \lcr (\mu-\nu-q+j)_i\rcr\prod_{k=0}^q \Gamma (\mu-\nu -k) (\mu-k-1)^{q-k},$$
et on conclut par le lemme \ref{Det1}.

\fin

\begin{Lemm} 
\label{Det3}
Pour tous entiers $n\ge q\ge 0$ et tous r\'eels $\mu,\nu$ tels que $\rho=\nu+n-\mu-2q>0,$ on a 
$${\det}^{}_q \lcr \lpa\frac{\Gamma(\rho+i+j)}{\Gamma(\rho +\mu+i+j)}\rpa(\nu-q+j)_{n-q+i}\rcr\; =\; F_q(\rho)\prod_{k=0}^q \lpa\frac{(\nu -k)_{n-q}}{\Gamma (\nu +n-k-q)}\rpa$$
avec $F_q(\rho) = {\det}^{}_q \lcr \Gamma(\rho+i+j)\rcr > 0.$ 
\end{Lemm}

\noindent
{\em Preuve} : En utilisant $\Gamma(t+1) = t \Gamma(t)$ on factorise
$${\det}^{}_q \lcr \lpa\frac{\Gamma(\rho+i+j)}{\Gamma(\rho +\mu+i+j)}\rpa(\nu-q+j)_{n-q+i}\rcr \; =\; \prod_{k=0}^q \lpa\frac{(\nu -k)_{n-q}}{\Gamma (\nu +n-k-q)}\rpa {\det}^{}_q \lcr \Gamma(\rho+i+j)\rcr.$$
Enfin, la stricte positivit\'e de ${\det}^{}_q \lcr \Gamma(\rho+i+j)\rcr$ d\'ecoule par exemple du th\'eor\`eme 5 de \cite{BI} et de la repr\'esentation int\'egrale  usuelle de la fonction $\Gamma$ sur $(0,+\infty),$ que l'on peut utiliser ici puisque $\rho > 0.$
\fin

\begin{Lemm} 
\label{Bloch}
Soit $A =\{ a_{ij}\}_{0\le i,j\le n}$ une matrice r\'eelle dont tous les coefficients sont non nuls. Pour tout $\delta\in(0,2n),$ soit la fonction $x\mapsto A^\delta (x) = \{a_{ij}\, x^{(\delta-i-j)_+}\}_{0\le i,j\le n}$ d\'efinie sur $\rl$ et soit $\rho_\delta = [\delta/2].$ Si les matrices
$$A_{11}^\delta \; =\; A\!\lcr\!\!\begin{array}{l} 
(1+\rho_\delta)\cdots\, n  \\
(1+\rho_\delta)\cdots\, n 
\end{array}\!\!\rcr\quad\mbox{et}\quad A_{22}^\delta \; =\; A\!\lcr\!\!\begin{array}{l}
1\cdots \rho_\delta \\
1\cdots \rho_\delta
\end{array}\!\!\rcr$$
sont inversibles, alors
$$\det A^\delta (x) \; \sim\; (\det A_{11}^\delta \det A_{22}^\delta)\, x^{(\rho_\delta+1)(\delta - \rho_\delta)}, \qquad x\to +\infty.$$
\end{Lemm}

\noindent
{\em Preuve} : Par hypoth\`ese on a $\rho_\delta < n$ et on peut \'ecrire par blocs
$$A^\delta (x)\; =\; \lcr\begin{matrix} A_{11}^\delta (x) & A_{12}^\delta (x)\\
A_{21}^\delta (x) & A_{22}^\delta (x)\end{matrix}\rcr$$
o\`u $A_{11}^\delta (x)$ est une matrice $(1+\rho_\delta)\times(1+\rho_\delta)$ et $A_{22}^\delta (x)$ une matrice $(n-\rho_\delta)\times(n-\rho_\delta).$ Le choix de la taille des blocs entra\^\i ne que $A_{22}^\delta (x) = A_{22}^\delta$ pour tout $x.$ Comme $A_{22}^\delta$ est inversible, on peut appliquer la formule classique
$$\det A^\delta (x) \; =\; \det A_{22}^\delta\times \det (A_{11}^\delta (x) - A_{12}^\delta (A_{22}^\delta)^{-1} A_{21}^\delta(x)).$$
On voit facilement que le coefficient $(i,j)$ de la matrice $A_{12}^\delta (A_{22}^\delta)^{-1} A_{21}^\delta(x)$ est un $O(x^{2(\delta -\rho_\delta -1) - i-j})$ quand $x\to +\infty$, quantit\'e n\'egligeable devant
$\lcr A_{11}^\delta (x)\rcr_{ij} =  a_{ij} x^{\delta - i - j}$
puisque $a_{ij}\neq 0$ et $\delta < 2(\rho_\delta +1).$ On en d\'eduit que
$$\det (A_{11}^\delta (x) - A_{12}^\delta (A_{22}^\delta)^{-1} A_{21}^\delta(x))\; \sim\; \det (A_{11}^\delta (x))\; =\;  (\det A_{11}^\delta)\, x^{(\rho_\delta+1)(\delta - \rho_\delta)}$$
par invertibilit\'e de $A_{11}^\delta,$ et o\`u l'\'egalit\'e r\'esulte de la formule de Leibniz. Tout ceci donne le comportement souhait\'e.

\fin

\section{Fin de la d\'emonstration du th\'eor\`eme}

Le cas ${\rm TP}_\infty$ ayant d\'ej\`a \'et\'e vu, on suppose dans toute cette partie que $(d,b)\notin\R$ et on montre le point (ii) du th\'eor\`eme, qui est le point d\'elicat. L'argument principal repose sur la notion de positivit\'e compl\`ete \'etendue. On commence par transformer le noyau ${\rm T}_{b,d}$ pour simplifier les calculs, en \'ecrivant 
$$e^{-xy}\Psi (b, b+1-d, xy)\; =\; y^{d-b} {\rm K}_{b,d}(x,y)$$ 
pour tout $x,y >0,$ avec
$${\rm K}_{b,d}(x,y)\; =\; \frac{1}{\Gamma (b)}\int_0^\infty e^{-x(t+y)} \frac{t^{b-1}}{(t +y)^{d}} dt.$$
Les fonctions d\'eterminantales
$$\Delta^n_{b,d}\; =\; {\det}^{}_n \lcr \frac{\partial^{i+j} {\rm K}_{b,d}}{\partial x^i \partial y^j}\rcr$$
d\'efinies sur $(0,+\infty)\times (0,+\infty),$ sont les quantit\'es centrales de notre probl\`eme. En effet, on sait par le th\'eor\`eme 2.2.6 p. 56 dans \cite{K} que pour tout $n\ge 1$ on a ${\rm K}_{b,d}\in {\rm RR}_{n+1}\Leftrightarrow{\rm T}_{b,d}\in {\rm TP}_{n+1}$ - voir \cite{K} p. 12 pour la notation RR - d\`es que
$$(-1)^{m(m+1)/2}\Delta^m_{b,d}\; >\; 0$$
 pour tout $m=1\ldots n.$ Mais comme $(d,b)\succ\E_n\,\Rightarrow\,(d,b)\succ\E_m$ pour tout $m=1\ldots n,$ on doit en fait consid\'erer \`a chaque fois un seul d\'eterminant. On voit en effet par r\'ecurrence imm\'ediate qu'il suffit de montrer
$$(d,b)\succ \E_n \; \Rightarrow\; (-1)^{n(n+1)/2}\Delta^n_{b,d}\; >\; 0$$
pour tout $n\ge 1$ afin d'avoir les inclusions d\'esir\'ees $(d,b)\succ \E_n \,\Rightarrow\,{\rm T}_{b,d}\in {\rm TP}_{n+1}.$  
 R\'eciproquement, on voit facilement par d\'efinition que si $\Delta^n_{b,d}$ change au moins une fois strictement de signe, alors on aura ${\rm K}_{b,d}\notin {\rm SR}_{n+1}\Leftrightarrow{\rm T}_{b,d}\notin {\rm SR}_{n+1}.$ 
 
On peut exprimer $(-1)^{n(n+1)/2}\Delta^n_{b,d}$ \`a l'aide d'un certain d\'eterminant de Hankel fonctionnel - un Turanien suivant la terminologie de \cite{KS} - avec poids. On voit d'abord par une suite d'op\'erations \'el\'ementaires sur les lignes utilisant les relations $\partial_x {\rm K}_{b,d}= - {\rm K}_{b, d-1}$ et $\partial_y {\rm K}_{b,d} = - x{\rm K}_{b, d}  - d{\rm K}_{b, d+1}$ que 
$$\Delta^n_{b,d}\; =\; {\det}^{}_n \lcr (d-j)_i\, {\rm K}_{b,d+i-j}(x,y)\rcr.$$
On en d\'eduit 
$${\rm sgn} ((-1)^{n(n+1)/2} \Delta^n_{b,d}(x,y))\; =\;{\rm sgn} ({\rm D}^n_{b,d}(xy))$$
pour tous $x,y >0,$ o\`u
$${\rm D}^n_{b,d}(z)\; = \; {\det}^{}_n \lcr (d-n+j)_i\, \Psi (b,b+1-d+n-i-j, z)\rcr$$
est bien un Turanien lest\'e par les $(d-n+j)_i.$ Signalons par ailleurs que sans ces poids, l'\'etude du signe serait imm\'ediate comme cons\'equence du th\'eor\`eme 2.8 dans \cite{IL}, lequel donne
$${\det}^{}_n \lcr \Psi (b,b+1-d+n-i-j, z)\rcr\; > \; 0$$
pour tout $b, z > 0, n\in\NN$ et $d\in\rl.$ Les complications calculatoires induites par la pr\'esence de poids dans les Turaniens sont \'evoqu\'ees dans l'introduction de \cite{KS}.\\

 La suite de cette partie est donc consacr\'ee \`a l'\'etude de la stricte positivit\'e de la fonction $z\mapsto{\rm D}^n_{b,d}(z)$ sur $(0,+\infty).$ Pour simplifier les \'ecritures, on utilisera parfois la notation
\begin{equation}
\label{fmunu}
f_{\mu,\nu} (z) \; =\; \Psi (\mu,\mu+1-\nu,z)\; =\; \frac{1}{\Gamma(\mu)}\int_0^\infty e^{-tz} \frac{t^{\mu -1}}{(1+t)^\nu} dt
\end{equation}
pour tout $\mu > 0$ et $\nu\in\rl.$

L'\'etude du signe de ${\rm D}^n_{b,d}$ est facile quand $d >n.$ En factorisant convenablement on voit en effet, avec la notation $h(a,c,x) = \Gamma (1+a -c) \Psi(a,c,x)$ pour tout $a > 0$ et $c< 1+a,$ que
\begin{eqnarray*}
{\rm D}^n_{b,d}(z) & = & \frac{1}{\prod_{k=0}^n \Gamma (d-k)}\times{\det}^{}_n \lcr h(b,b+1-d-n+i+j, z)\rcr.
\end{eqnarray*}
La formule exacte donn\'ee dans \cite{BI} p. 12 montre alors que $z\mapsto {\det}^{}_n \lcr h(b,b+1-d-n+i+j, z)\rcr$ est strictement positive (et m\^eme compl\`etement monotone) sur $(0,+\infty)$ puisque $d > n,$ et il en est donc de m\^eme pour ${\rm D}^n_{b,d}(z).$ Ceci entra\^\i ne donc la propri\'et\'e ${\rm TP}_{n+1}$ pour ${\rm T}_{b,d}$ lorsque $d >n,$ ce que nous avions d\'ej\`a vu dans la partie pr\'ec\'edente par un argument de sym\'etrie.

L'\'etude devient cependant plus ardue quand $d < n,$ puisque certaines fonctions $h$ dans le d\'eterminant ont alors des param\`etres n\'egatifs. On peut certes reprendre formellement les calculs du th\'eor\`eme 5 de \cite{BI} \`a l'aide de la transformation de Kummer et d'une repr\'esentation int\'egrale de $\Psi$ le long d'un chemin de type Hankel - voir la formule 6.11.2 (9) p.273 dans \cite{E}, mais ceci exprime ${\det}^{}_n \lcr h(b,b+1-d-n+i+j, z)\rcr$ comme une int\'egrale curviligne multiple dont il para\^\i t malais\'e de d\'eterminer le signe. D'une fa\c{c}on g\'en\'erale, il ne semble pas que l'on dispose de transformations d\'eterminantales simples, telle celle donn\'ee dans l'exercice II.1.68 de \cite{PS}, qui fasse appara\^\i tre naturellement la barri\`ere $\E_n.$  

\subsection{Le cas $n=1$} On traite ce cas \`a part pour simplifier l'exposition de la m\'ethode, qui sera la m\^eme dans le cas g\'en\'eral avec des arguments plus \'elabor\'es. Par le r\'esultat de Dette et Munk, on sait que ${\rm T}_{b,d}\in{\rm TP}_2$ si $b > 1$ et on est donc ramen\'e \`a une \'etude sur la bande $\{0< b<1, d\in\rl\}.$ On r\'eduit encore le domaine d'\'etude en consid\'erant les carr\'es  
$$\E_1^+ = \{0<x,y<1\}\quad\mbox{et}\quad \E_1^- = \{-1\le x <0, 0<y<1\}.$$
Ces carr\'es sont ouverts sauf $\E_1^-$ \`a son ar\^ete gauche $\{x =-1, 0<y<1\}.$ L'identit\'e en loi (\ref{crescent}) et le caract\`ere ${\rm TP}_\infty$ du noyau multiplicatif associ\'e \`a $\U$ montrent automatiquement  que si ${\rm T}_{b,d}$ est ${\rm TP}_2$ pour $(d,b)\in\E_1^+$ alors il le sera aussi sur $\{0< b<1, d>0\},$ et que si ${\rm T}_{b,d}$ n'est pas ${\rm SR}_2$ pour $(d,b)\in\E_1^-$ alors il ne le sera pas non plus sur $\{0< b<1, d<0\}.$ \\

Il y a deux fa\c{c}ons de montrer que ${\rm T}_{b,d}$ est ${\rm TP}_2$ sur $\E_1^+.$ La premi\`ere, propre au cas $n=1$, utilise la notion de compl\`ete monotonicit\'e hyperbolique (HCM) faisant l'objet du livre \cite{Bo} - voir en particulier le chapitre 4.3 et le chapitre 5. En utilisant la d\'efinition de ${\rm T}_{b,d}$, la caract\'erisation de ${\rm TP}_2$ par log-concavit\'e et la log-concavit\'e de $x\mapsto e^{-e^x}$ sur $\rl,$ il suffit en effet de d\'emontrer que $x\mapsto f_{b,d}(e^x)$ est log-concave sur $\rl,$ ce qui est encore \'equivalent  - voir le paragraphe 6.4 pp. 101-102 dans \cite{Bo} - \`a la propri\'et\'e HM pour la fonction $x\mapsto f_{b,d}(x)$ sur $(0,+\infty).$ On voit cependant par (\ref{fmunu}) et la propri\'et\'e 5.1(vii) p. 68 dans \cite{Bo} que cette derni\`ere est HCM (et donc HM) pour $d>0$ puisque la fonction $t\mapsto t^{b-1}(1+t)^{-d}$ est alors HCM - voir la propri\'et\'e 5.1(ii) p. 68 dans \cite{Bo}.

La deuxi\`eme mani\`ere montre la stricte positivit\'e de ${\rm D}^1_{b,d}$ pour tout $(d,b)\in\E_1^+,$ et c'est cette argumentation qui sera reprise pour $n\ge 2.$ On remarque d'abord par l'in\'egalit\'e de H\"older que la fonction $\nu \mapsto f_{b,\nu} (z)$ est log-convexe sur $\rl$ pour tout $z > 0.$ On \'ecrit ensuite
 \begin{eqnarray*}
{\rm D}^1_{b,d} & = & \lva \begin{matrix} 
f_{b,d-1}  & f_{b,d} \\
(d-1)f_{b,d} & d f_{b,d+1}\end{matrix}\rva\; =\; d(f_{b,d-1}f_{b,d+1} - f^2_{b,d}) \, +\, f^2_{b,d}
\end{eqnarray*}
et cette fonction est bien strictement positive sur $(0,+\infty)$ puisque $d >0$ et $f_{b,d}$ ne s'annule pas. 

\begin{Rq} {\em Dans chacune des deux preuves on n'a pas utilis\'e l'hypoth\`ese $(d,b)\in\E_1^+$ mais seulement $d >0.$ Dans la suite, on aura cependant besoin des domaines $\E_n^+$ correspondant \`a $n\ge 2.$}
\end{Rq}
 
On montre maintenant que ${\rm D}^1_{b,d}$ change strictement de signe sur $(0,+\infty)$ pour $(d,b)\in\E_1^-.$ Il est imm\'ediat que $f_{b,\nu} (z)\sim z^{-b}$ quand $z\to +\infty$ pour tout $\nu\in\rl,$ ce qui entra\^\i ne
$${\rm D}^1_{b,d}(z) \; \sim \; z^{-2b}$$ 
en $+\infty$ et montre donc que ${\rm D}^1_{b,d}(z) > 0$ pour $z$ assez grand. On \'etablit que ${\rm D}^1_{b,d}(z)\to -\infty$ quand $z\to 0+$ \`a l'aide du lemme suivant, \'el\'ementaire, et qui jouera \'egalement un r\^ole dans le cas g\'en\'eral $n\ge 2.$

\begin{Lemm} \label{Asymph} Pour tout $\mu > 0$ et $\nu\in\rl$ on a 
$$f_{\mu,\nu} (z)\; \sim\; \lacc\begin{array}{ll}
{\displaystyle \frac{\Gamma(\nu -\mu)}{\Gamma (\nu)}} & \mbox{si $\mu < \nu$} \\
{\displaystyle \frac{\Gamma(\mu -\nu)}{\Gamma (\mu)} z^{\nu-\mu}}& \mbox{si $\mu > \nu$} \\
{\displaystyle \frac{1}{\Gamma (\mu)} (-\log z)} & \mbox{si $\mu = \nu$} 
\end{array}\right.$$
quand $z\to 0\!+\!.$ 
\end{Lemm}

On d\'ecoupe $\E_1^-$ en un triangle sup\'erieur $\E_1^{s-} = \{-1\le x < y -1<0\},$ un triangle inf\'erieur $\E_1^{i-} = \{-1<y-1 < x<0\},$ et une diagonale $\E_1^{d-} = \{-1<x = y -1<0\}.$ Pour $(d,b)\in\E_1^{s-}$ on trouve \`a l'aide du lemme \ref{Asymph} et apr\`es simplification le comportement asymptotique
$${\rm D}^1_{b,d}(z)\; \sim\; (b-1)(b-d-1)\lpa\frac{\Gamma(b-d-1)}{\Gamma(b)}\rpa^2 z^{2(d-b)}\; \to\; -\infty$$
quand $z\to 0\!+\!.$ Quand $(d,b)\in\E_1^{i-}$ et $(d,b)\in\E_1^{d-}$ ces comportements deviennent respectivement 
$$\frac{d\,\Gamma(b-d+1)\Gamma(d+1-b)}{\Gamma(b)\Gamma(d+1)}\, z^{d-1-b}\; \to\; -\infty\quad\mbox{et}\quad \frac{d\,\Gamma(b-d+1)}{\Gamma(b)^2}\, (-\log z)z^{-2}\; \to\; -\infty.$$ 
Tout ceci prouve que ${\rm D}^1_{b,d}$ change strictement de signe sur $(0,+\infty),$ et donc que ${\rm T}_{b,d}$ n'est pas ${\rm SR}_2,$ quand $(d,b)\in\E_1^-.$

\subsection{Le cas g\'en\'eral} On fixe $n\ge 2$ et on introduit dans $\HH$ les deux escaliers adjacents suivants, constitu\'es chacun de $n$ carr\'es ouverts disjoints :
$$\E_n^-\,=\,\bigcup_{k=0}^{n-1}\{k-1< x< k, n-k-1<y<n-k\}\;\;\;\mbox{et}\;\;\;\E_n^+\,=\,\bigcup_{k=1}^{n}\{k-1< x< k, n-k<y<n-k+1\}.$$
 On sait d\'ej\`a que ${\rm T}_{b,d}$ est ${\rm TP}_{n+1}$ si $b >n.$ A nouveau, l'identit\'e en loi (\ref{crescent}) et le caract\`ere ${\rm TP}_\infty$ du noyau multiplicatif associ\'e \`a $\U$ montrent que si ${\rm T}_{b,d}$ est ${\rm TP}_{n+1}$ pour $(d,b)\in\E_n^+$ alors il le sera aussi pour $(d,b)\succ \E_n$ et $b\le n.$ De m\^eme, si ${\rm T}_{b,d}$ n'est pas ${\rm SR}_{n+1}$ pour $(d,b)\in\E_n^-$ alors il ne le sera pas non plus pour $(d,b)\prec\E_n$ et $-d\notin\NN^*.$\\

On montre d'abord par r\'ecurrence que pour $(d,b)\in \E_n^+,$ on a ${\rm D}^n_{b,d} > 0$ sur $(0,+\infty).$ Rappelons que le cas $n=1$ a d\'ej\`a \'et\'e vu.\\

\noindent
{\bf Cas $n=2.$} On isole ce cas \`a nouveau pour simplifier l'exposition. On a 
$${\rm D}^2_{b,d} \; = \; \lva \begin{matrix} 
f_{b,d-2}  & f_{b,d-1} & f_{b,d}\\
(d-2)f_{b,d-1} &(d-1)f_{b,d} & d f_{b,d+1}\\
(d-2)(d-1)f_{b,d} &d(d-1)f_{b,d+1} & d(d+1) f_{b,d+2}\end{matrix}\rva$$
et $\E_2^+ = \{0< b <1<d<2\}\cup\{0< d <1<b<2\}.$ Par un argument de log-convexit\'e semblable au cas $n=1$, il est possible de montrer directement la stricte positivit\'e de ${\rm D}^2_{b,d}$ pour $\{0< d <1<b<2\}$ en d\'eveloppant par rapport \`a la troisi\`eme ligne, et pour $\{0< b <1<d<2\}$ en d\'eveloppant par rapport \`a la premi\`ere ligne. Mais cet argument ne se g\'en\'eralise pas aux plus grandes valeurs de $n.$ On utilise plut\^ot l'identit\'e d\'eterminantale de Sylvester - voir le chapitre 0 p. 3 dans \cite{K} - avec pour bloc pivot la matrice scalaire $(d-1)f_{b,d}$. On obtient
$$(d-1)f_{b,d} {\rm D}^2_{b,d} \; = \; \lva \begin{matrix} 
d(d-1) {\rm D}^1_{b,d+1} & (d-1)(d-2) {\rm D}^1_{b,d}\\
{\rm D}^1_{b,d} & {\rm D}^1_{b,d-1}
\end{matrix}\rva,$$
ce qui entra\^ine
$${\rm D}^2_{b,d} \; = \; \frac{1}{f_{b,d}} (d\, {\rm D}^1_{b,d-1}
{\rm D}^1_{b,d+1} + (2-d) ({\rm D}^1_{b,d})^2)$$
qui est bien strictement positive pour $(d,b) \in\L^+_2$ par l'hypoth\`ese de r\'ecurrence, puisque l'on sait que ${\rm D}^1_{b,d-1}{\rm D}^1_{b,d+1} \ge 0, f_{b,d}> 0$ et ${\rm D}^1_{b,d} >0$ pour $(d,b)\in\E_2^+.$  \\

\noindent
{\bf Cas $n\ge 3.$} On d\'eveloppe ${\rm D}^n_{b,d}$ de nouveau par l'identit\'e de Sylvester, en choisissant cette fois pour bloc pivot la sous-matrice centrale de taille $n-1$, qui est donc bord\'ee par les lignes et les colonnes num\'ero 0 et $n.$ Apr\`es simplifications, il vient 
$${\rm D}^n_{b,d} \; = \; \frac{1}{{\rm D}^{n-2}_{b,d}} (d\, {\rm D}^{n-1}_{b,d-1}{\rm D}^{n-1}_{b,d+1} + (n-d) ({\rm D}^{n-1}_{b,d})^2).$$
Ceci permet de conclure par r\'ecurrence, puisque l'on sait que ${\rm D}^{n-1}_{b,d-1}{\rm D}^{n-1}_{b,d+1} \ge 0, {\rm D}^{n-2}_{b,d} > 0$ et ${\rm D}^{n-1}_{b,d} > 0$ pour $(d,b)\in\E_n^+.$ \\

Tout ceci montre l'inclusion recherch\'ee $(d,b)\succ \E_n\;\Rightarrow\; {\rm T}_{b,d}\in {\rm TP}_{n+1}$ pour tout $n\ge 1.$

\begin{Rq}{\em L'argument pr\'ec\'edent \'etant uniquement alg\'ebrique, il est tentant de chercher une preuve probabiliste du sens direct de la partie (ii) du th\'eor\`eme par le biais d'identit\'es en loi multiplicatives. Si on divise $\E_{n+1}^-$ en $(n+1)$ marches disjointes $\E_{n+1}^{\nu-}$ num\'erot\'ees dans l'ordre descendant de $\nu = 0$ \`a $n,$ on a $\E_n^+ = \E_{n+1}^{1-}\cup\ldots\cup\E_{n+1}^{n-}$ et on sait que ${\rm T}_{b,d}\in {\rm TP}_{n+1}$ pour $(d,b)\in \E_{n+1}^{0-}\subset\{ b>n\}.$ D'autre part, la descente de $\E_{n+1}^{\nu-}$ \`a $\E_{n+1}^{\nu+1-}$ correspond \`a une factorisation multiplicative pour tout $\nu =0, \ldots, n-1$. En effet, on a 
$$\EE[\T(a,b,c)^s]\; =\; \frac{c(c+1)(a+b+s)}{(a+b)(c+s)(c+1+ s)}\times \EE[\T(a,b+1,c+2)^s]$$
pour tout $s\in\CC$ tel que $\Re(s) > (-a)\vee(-c),$ 
et le quotient \`a droite est la transform\'ee de Mellin d'une certaine variable al\'eatoire positive si $c < a+b$ autrement dit $d > 0.$ De fa\c{c}on plus pr\'ecise, si on pose $\U_\rho = \eps_\rho + (1-\eps_\rho) \U$ o\`u $\rho = c(a+b)^{-1} < 1$ et $\eps_\rho$ est une variable ind\'ependante de Bernoulli de param\`etre $\rho$, on trouve la factorisation
$$\T(a,b,c)\;\elaw\;\U_\rho^{1/c}\,\times\, \U^{1/(c+1)}\, \times\, \T(a,b+1,c+2)$$
pour tout $a,b,c > 0$ tels que $c<a+b.$ Cette identit\'e permet bien pour tout $\nu =0, \ldots, n-1$ de passer de $\E_{n+1}^{\nu-}$ \`a $\E_{n+1}^{\nu+1-}$ par convolution multiplicative, en choisissant convenablement les param\`etres. Si le noyau multiplicatif associ\'e \`a $\U_\rho^{1/c}\,\times\, \U^{1/(c+1)}$ avait une positivit\'e compl\`ete appropri\'ee, on pourrait esp\'erer d\'eduire la propri\'et\'e ${\rm TP}_{n+1}$ pour ${\rm T}_{b,d}$ sur tout $\E_n^+$ de celle d\'ej\`a connue sur $\E_{n+1}^{0-}$ en faisant varier $\nu,$ descendant l'escalier. Malheureusement, la densit\'e de $\U_\rho^{1/c}\,\times\, \U^{1/(c+1)}$ est la fonction
$$\rho (c+1) x^{c-1} (d + (1-d)x)\Un_{(0,1)}(x),$$
de sorte que le noyau multiplicatif associ\'e n'est pas ${\rm SR}_2$ si $d< 1$ et n'est pas ${\rm SR}_3$ si $d > 1.$}
\end{Rq}

On \'etablit maintenant que ${\rm D}^n_{b,d}$ change au moins une fois strictement de signe sur $(0,+\infty)$ quand $(d,b)\in\L^-_n.$ Par la discussion pr\'ec\'edente, ceci montrera l'inclusion 
$$(d,b) \prec \E_n\;\mbox{et}\;-\!d\notin\NN^*\;\Rightarrow\; {\rm T}_{b,d}\notin {\rm SR}_{n+1}.$$ 
Les arguments ressemblent \`a ceux du cas $n=1$ mais sont h\'elas plus techniques. Par le lemme \ref{Det1}, on voit que 
\begin{equation}
\label{Asu}
{\rm D}^n_{b,d}(z)\;\sim\; (\prod_{k=0}^n k!)z^{-(n+1)b}
\end{equation}
quand $z\to +\infty,$ de sorte que pour tout $b >0, d\in\rl$ on a ${\rm D}^n_{b,d}(z) > 0$ si $z$ est assez grand. Il suffit donc de montrer  
$$\lim_{z\to0+}{\rm D}^n_{b,d}(z)\; \to\; -\infty$$ 
pour tout $(d,b)\in \E_n^-.$ Dans le demi-plan $\{ b>0\},$ l'ensemble $\E_n^-$ est constitu\'e des $n$ carr\'es ouverts disjoints
$$\E_{n,p}^-\; =\;\{p-2<x<p-1\}\,\cap\, \{n-p<y<n+1-p\}$$ 
pour $p= 1, \ldots, n,$ qu'on divise en un triangle sup\'erieur $\E_{n,p}^{s,-} = \E_{n,p}^-\cap\{x < y +2(p-1) -n\},$ un triangle inf\'erieur $\E_{n,p}^{i,-} = \E_{n,p}^-\cap\{x > y +2(p-1) -n\},$ et une diagonale $\E_{n,p}^{d,-} = \E_{n,p}^-\cap\{x = y +2(p-1) -n\}.$ \\

Pour $(d,b)\in\E_{n,p}^{s,-}$, le lemme \ref{Asymph} entra\^\i ne avec la notation du lemme \ref{Bloch} que
$${\rm D}^n_{b,d}(z)\; \sim\; \det A^\delta (1/z)$$
lorsque $z\to 0+,$ o\`u l'on a pos\'e $\delta = b+n-d > 1$ et o\`u les coefficients de la matrice $A$ sous-jacente sont donn\'es par
$$a^{}_{ij}\; =\; \lacc\begin{array}{ll} {\displaystyle (d-n+j)_i \frac{\Gamma(b+n-d-i-j)}{\Gamma(b)}} & \mbox{si $i+j\le 2(n+1-p),$}\\
 {\displaystyle (d-n+j)_i \frac{\Gamma(d-n-b+i+j)}{\Gamma(d-n+i+j)}} & \mbox{si $i+j > 2(n+1-p).$}\end{array}\right.$$
Avec la notation du lemme \ref{Bloch}, on voit que $\rho_\delta = n+1-p < \delta.$ En posant $q = n+1-p, \mu = b$ et $\nu = d+1 -p,$ on a $\mu > \nu +q$ par d\'efinition de $\E_{n,p}^{s,-}.$ Par le lemme \ref{Det2}, le d\'eterminant de la matrice $A^\delta_{11}$ sous-jacente vaut donc
$${\det}^{}_q A^\delta_{11}\; =\; \prod_{k=0}^q k!\,\Gamma (\mu-\nu -k) (\mu-k-1)^{q-k}\; <\; 0,$$
o\`u l'in\'egalit\'e stricte vient de $q-1<\mu<q.$ Si $p = 1$ on a $A^\delta = A^\delta_{11}$ et la formule de Leibniz donne $\det A^\delta (1/z) = \det A^\delta_{11} z^{-(n+1)\delta} \to  -\infty$ quand $z\to 0\!+\!.$ Si $p > 1,$ les coefficients de la matrice $A^\delta_{22}$ sous-jacente qui est de taille $r=p-2,$ s'\'ecrivent
$$(d-r+j)_{n-r+i} \frac{\Gamma(d+n-b-2r+i+j)}{\Gamma(d+n-2r+i+j)}$$ pour $0\le i,j\le r.$ En posant $\rho = d+n -b-2r,$ on a $\rho > 0$ par d\'efinition de $\E_{n,p}^{-}.$ Le lemme \ref{Det3} donne 
$$\det A^\delta_{22}\; =\; F_q(\rho)\prod_{k=0}^r \lpa\frac{(d -k)_{n-r}}{\Gamma (d +n-k-r)}\rpa\; >\; 0,$$
o\`u l'in\'egalit\'e stricte vient de $d > p-2.$ On a d'autre part $\rho_\delta < \delta$ et le lemme \ref{Bloch} entra\^\i ne finalement 
$${\rm D}^n_{b,d}(z) \; \sim\;  (\det A_{11}^\delta \det A_{22}^\delta) z^{(\rho_\delta+1)(\rho_\delta-\delta)}\;\to\; -\infty$$
quand $z\to 0\!+\!.$\\

L'argumentation est la m\^eme sur $\E_{n,p}^{i,-},$ sauf que l'on voit par le lemme \ref{Asymph} que la matrice $A_{22}^\delta$ sous-jacente est cette fois-ci toujours d\'efinie, et de taille $p-1$. On calcule
$${\det}^{}_{n-p} (A^\delta_{11})\; =\; \prod_{k=0}^{n-p} k!\,\Gamma (b-d-1+p -k) (b-k-1)^{n-p-k}\; >\; 0$$
o\`u l'in\'egalit\'e stricte vient de $b>n-p,$ et
$${\det}^{}_{p-1} (A^\delta_{22})\; =\; F_{p-1}(d+n -b+2-2p)\prod_{k=0}^{p-1} \lpa\frac{(d -k)_{n+1-p}}{\Gamma (d +n-k+1-p)}\rpa\; <\; 0,$$
o\`u l'in\'egalit\'e stricte vient de $p-2<d<p-1$ et de $n+1-p\ge 1.$ On a toujours $\rho_\delta < \delta$ et on en d\'eduit ${\rm D}^n_{b,d}(z)\to -\infty$ quand $z\to 0+$ par le lemme \ref{Bloch}.\\

Sur la diagonale $\E_{n,p}^{d,-},$ la preuve est l\'eg\`erement diff\'erente et repose sur une variante du lemme \ref{Bloch} que nous laissons au lecteur. Pour $p = 1$ on trouve \`a l'aide des lemmes \ref{Det2} et \ref{Asymph}
$${\rm D}^n_{b,d}(z)\; \sim\; \lpa\prod_{k=0}^{n-1} k!\,(d+k)\Gamma (n-k) 
(d+n-k-1)^{n-1-k}\rpa (-\log z)\, z^{-n^2}\; \to\; -\infty$$
quand $z\to 0+,$ o\`u la stricte n\'egativit\'e du coefficient directeur vient de $-1<d<0.$  Pour $p\ge 2$ on trouve \`a l'aide des lemmes \ref{Det2}, \ref{Det3} et \ref{Asymph} 
\begin{eqnarray*}
{\rm D}^n_{b,d}(z)& \sim & \lpa\prod_{k=0}^{n+1-p} (d+1-p+k)\rpa \lpa\prod_{k=0}^{n-p} k!\,\Gamma (b-d-1+p -k) (b-k-1)^{n-p-k}\rpa\; \times\\
& & F_{p-2}(d+n -b-2p+4) \lpa\prod_{k=0}^{p-2}\frac{(d -k)_{n+2-p}}{\Gamma (d +n-k+2-p)}\rpa (-\log z)\, z^{-(n+1-p)^2}\; \to\; -\infty
\end{eqnarray*}
quand $z\to 0+,$ o\`u la stricte n\'egativit\'e du coefficient directeur est celle du premier facteur et vient de $p-2<d<p-1.$\\

Pour terminer la preuve de la partie (ii) du th\'eor\`eme, il reste \`a montrer  l'inclusion 
$$(d,b) \prec \E_n\;\mbox{et}\;-\!d\in\NN^*\;\Rightarrow\; {\rm T}_{b,d}\notin {\rm SR}_{n+1}.$$ 
Par (\ref{crs}), (\ref{crescent}) et le raisonnement pr\'ec\'edent, ceci est une cons\'equence de
$$\lim_{z\to 0+}{\rm D}^n_{b,d}(z)\; =\; -\infty$$ 
pour $d=-1$ et $n-1<b<n.$ Mais ce dernier point a en fait d\'ej\`a \'et\'e d\'emontr\'e dans le sous-cas $(d,b)\in\E_{n,1}^{s,-}$ de la discussion pr\'ec\'edente.

\fin

\begin{Rqs} {\em (a) Par des arguments semblables, on montre que si $(d,b)\in \E_n^+,$ alors ${\rm D}^n_{b,d}(z) \to +\infty$ quand $z\to 0\!+\!.$

(b) Dans le cas $d > 0,$ le th\'eor\`eme fournit un exemple de noyaux n'appartenant pas \`a ${\rm SR}_\infty$ et qui sont obtenus par multiplication ponctuelle de deux noyaux ${\rm TP}_\infty.$ En effet, le noyau $e^{xy^{-1}}$ est ${\rm TP}_\infty$ sur $(0,+\infty)\times (0, +\infty)$ et par convolution multiplicative le noyau $\Psi(b, b+1-d, xy^{-1})$ est aussi ${\rm TP}_\infty$ quand $d >0,$ puisque la fonction $t\mapsto (1+e^t)^{-d}$ est alors ${\rm PF}_\infty$ sur $\rl$ - voir pour cela le th\'eor\`eme 3.2 dans \cite{DM}. En revanche, on a d\'emontr\'e que le noyau
$${\rm T}_{b,d}(xy^{-1})\; =\; e^{xy^{-1}}\!\times\Psi(b, b+1-d, xy^{-1})$$
n'est pas ${\rm SR}_\infty$ sur $(0,+\infty)\times (0, +\infty)$ si $(d,b)\notin\R.$ Dans cet exemple, la seule propri\'et\'e qui reste stable par multiplication ponctuelle est la propri\'et\'e ${\rm TP}_2.$

(c) Le th\'eor\`eme montre que le noyau multiplicatif de $\T(a,b,c)$ n'est pas ${\rm SR}_{[b+2]}$ si $b\notin\NN^*$ et $c > a+b.$ Ceci entra\^\i ne que le noyau multiplicatif de $\B(a,b)$ n'est pas ${\rm SR}_{[b+2]}$ si $b\notin\NN^*$ puisque celui de $\G(c)$ est toujours ${\rm TP}_\infty.$ Il est aussi possible de voir directement cette propri\'et\'e \`a l'aide du noyau $(y-x)_+^{b-1},$ t\^ache que nous laissons au lecteur.}
\end{Rqs}

\section{Autres d\'eterminants associ\'es \`a $\Psi.$}

Dans cette derni\`ere partie on donne quelques remarques sur des Wronskiens et des Turaniens associ\'es \`a la fonction de Tricomi, dans l'esprit de \cite{KS}. 

\subsection{Wronskiens} On pose 
$$z\; \mapsto\; {\rm W}^n_{b,d} (z) \; =\; {\det}^{}_n\lcr \frac{\partial^i}{\partial z^i}\Psi(d-j,d+1-b-j, z)\rcr$$
sur $(0,+\infty).$ On voit par la formule de Leibniz et les formules 6.5(6) et 6.6(12) pp.257-258 dans \cite{E} que 
\begin{equation}
\label{Asd}
{\rm D}^n_{b,d} (z)\; =\; z^{(n+1)(d-b)} {\rm W}^n_{b,d} (z).
\end{equation}
Ainsi, le th\'eor\`eme entra\^\i ne entre autres que
$${\rm W}^n_{b,d} > 0\;\Leftrightarrow\; (d,b)\,\succ\,\E_n \;\mbox{ou}\; (d,b)\in\R$$
pour tout $n\ge 1.$ Il ne semble pas possible d'obtenir une telle caract\'erisation de la constance du signe du Wronskien ${\rm W}^n_{b,d}$ \`a l'aide de l'\'equation hyperg\'eom\'etrique. La transformation
$$\Psi (b, b+1-d, x) \; = \; e^{x/2} x^{-1/2-\mu} W_{\kappa, \mu} (x)$$
o\`u $\kappa = (d-1-b)/2$ et $\mu = (b-d)/2$ - voir 6.8(4) p. 264 dans \cite{E} - et l'\'equation diff\'erentielle de type Sturm-Liouville v\'erifi\'ee par la fonction de Whittaker $W_{\kappa, \mu}$ ne semblent pas non plus pouvoir aider. 

Cette \'ecriture wronskienne permet de pr\'eciser le comportement de $z\mapsto {\rm D}^n_{b,d} (z)$ pour $(d,b) \succ \E_n$ et $b\ge d > 0.$ On voit en d\'erivant colonne par colonne et en utilisant \`a nouveau la formule de Leibniz et les formules 6.5(6) et 6.6(12) pp. 257-258 dans \cite{E} que
$$\frac{{\rm d}}{{\rm dz} }{\rm W}^n_{b,d} (z)\; =\; -d z^{(n+1)(b-d)+1} {\rm {\tilde D}}^n_{b,d} (z)$$
o\`u ${\rm {\tilde D}}^n_{b,d} (z)$ est obtenu \`a partir de ${\rm D}^n_{b,d} (z)$ en rempla\c{c}ant $d$ par $d+1$ dans toute la derni\`ere colonne. Un argument en tous points semblable \`a celui du sens direct de la partie (ii) du th\'eor\`eme montre alors que ${\rm {\tilde D}}^n_{b,d} > 0$ pour $(d,b) \succ \E_n.$ On en d\'eduit que $z\mapsto {\rm W}^n_{b,d} (z)$ est strictement d\'ecroissante pour $(d,b) \succ \E_n$ et $d >0,$ et strictement croissante pour $(d,b) \succ \E_n$ et $d <0.$ Ces propri\'et\'es qui s'accordent avec (\ref{Asu}) et (\ref{Asd}), sont bien s\^ur un peu attendues pour un Wronskien. On en d\'eduit aussi $z\mapsto {\rm D}^n_{b,d} (z)$ est strictement d\'ecroissante pour $(d,b) \succ \E_n$ et $b\ge d >0,$ en accord avec la derni\`ere remarque (a) ci-dessus. 

Il est plausible que $z\mapsto {\rm D}^n_{b,d} (z)$ est strictement d\'ecroissante pour $(d,b) \succ \E_n$ et $d >0.$ Rappelons que si $d > n,$ la transformation de Kummer et le th\'eor\`eme 5 dans \cite{BI} entra\^\i nent que cette fonction est compl\`etement monotone.

\subsection{Turaniens.} Comme leur nom l'indique, les Turaniens servent \`a montrer des in\'egalit\'es de type Tur\'an. Introduites apr\`es-guerre pour les polyn\^omes de Legendre, ces derni\`eres peuvent \^etre vues comme des variations sur la log-convexit\'e ou la log-concavit\'e de fonctions sp\'eciales en le param\`etre. Une \'etude syst\'ematique pour de nombreuses classes de polyn\^omes orthogonaux a \'et\'e men\'ee dans \cite{KS}. L'int\'er\^et port\'e \`a ces in\'egalit\'es pour des fonctions non polynomiales est en revanche plus r\'ecent et nous renvoyons entre autres aux articles \cite{B,BI,IL} ainsi qu'\`a leur bibliographie. 


Certaines in\'egalit\'es de Tur\'an sont tr\`es faciles. Ainsi, il est imm\'ediat par l'in\'egalit\'e de H\"older et d\'ej\`a mentionn\'e plus haut que 
$$\Psi^2(a,c)\; \le \; \Psi(a,c-1)\Psi(a,c+1)$$
sur $(0,+\infty)$ pour tout $a > 0$ et $c\in\rl.$ Cependant, on sait que $\Psi^2(a,c)\sim\Psi(a,c-1)\Psi(a,c+1)$ en $+\infty$, et on peut donc se poser la question de la taille de l'\'ecart dans l'in\'egalit\'e. A cet effet, le th\'eor\`eme 4 dans \cite{BI} montre que
$$\Psi(a,c-1)\Psi(a,c+1)\, - \,\Psi^2(a,c)\; \le \; \frac{a}{c(c-a-1)}\,\Psi^2(a,c)$$
sur $(0,+\infty)$ pour tout $a > 0> c$ et que
$$\Psi(a,c-1)\Psi(a,c+1)\, - \,\Psi^2(a,c)\; \le \; \frac{1}{c-2}\,\Psi^2(a,c)$$
sur $(0,+\infty)$ pour tout $a > c-1 > 1.$ L'argument repose essentiellement sur une repr\'esentation int\'egrale pour le quotient de fonctions de Tricomi ayant des param\`etres contigus, d\^ue \`a Ismail et Kelker - voir (3.6) dans \cite{BI}. Il est possible de compl\'eter ces deux in\'egalit\'es \`a l'aide de notre r\'esultat principal.

\begin{Corru} 
Supposons $a \ge 1$ ou bien $a > 0$ et $c\le a+2.$ On a
\begin{equation}
\label{Turan}
\Psi(a,c-1)\Psi(a,c+1)\, - \,\Psi^2(a,c)\; \le \; \frac{1}{x}\,\Psi(a,c)\Psi(a,c-1)
\end{equation}
sur $(0,+\infty).$ 
\end{Corru}

\noindent
{\em Preuve} : En posant $d=2+a-c,$ on a $\{a \ge 1\}\cup\{a > 0, c\le a+2\} = \{(d,a)\succ \L_1^+\}.$ Par le th\'eor\`eme, on sait que la fonction $g_{a,c}(e^x)$ o\`u $g_{a,c}(x) = e^{-x}\Psi(a, c-1,x)$ est alors ${\rm PF}_2$ sur $\rl,$ autrement dit 
$$g_{a,c}(x^2g_{a,c}''+xg_{a,c}')\; \le \; (xg_{a,c}')^2$$
sur $(0,+\infty).$ On obtient alors ais\'ement (\ref{Turan}) en utilisant $g_{a,c}' = -g_{a,c+1}$ et $g_{a,c}'' = g_{a,c+2}.$
\fin

\begin{Rqs} {\em (a) Le domaine de validit\'e de (\ref{Turan}) contient strictement celui du th\'eor\`eme 4 de \cite{BI}, et autorise par exemple des valeurs de $c$ dans $[0,2].$ D'autre part l'in\'egalit\'e (\ref{Turan}) est meilleure que les deux pr\'ec\'edentes sur leurs domaines respectifs, pour peu que $x$ soit suffisamment \'eloign\'e de z\'ero. En effet, il est \'evident par d\'efinition de $\Psi$ que $c(c-a-1)\Psi(a,c-1,x) <  ax \Psi(a,c,x)$
pour tout $a> 0 > c$ et $x > c(c-a-1)/a,$ et que 
$(c-2)\Psi(a,c-1,x)<  x\Psi(a,c,x)$
pour tout $a> c-1>1$ et $x > c-2.$ Signalons enfin, toujours par notre r\'esultat principal, que (\ref{Turan}) est fausse si $a < 1$ et $c > a+1.$ 

(b) Dans la preuve de ce corollaire on a utilis\'e de mani\`ere non travaill\'ee le caract\`ere ${\rm TP}_2$ du noyau ${\rm T}_{b,d},$ sans la simplification par ${\rm K}_{b,d}$ qui permet d'enlever la variable $x$ dans le d\'eterminant et de d\'emontrer la partie (ii) du th\'eor\`eme. Il est possible de faire la m\^eme chose pour des d\'eterminants d'ordre sup\'erieur. On obtient par exemple par ${\rm TP}_3$ des in\'egalit\'es h\'elas pas tr\`es parlantes, avec produits triples et param\`etres $c$ contigus jusqu'\`a l'ordre 2.}
\end{Rqs}

Si on suppose maintenant $a > 1$ et $c < 1+a,$ la transformation de Kummer et l'in\'egalit\'e $\Psi^2(1+a-c, 2-c)\le \Psi(1+a-c, 1-c)\Psi(1+a-c, 3-c)$ entra\^\i nent  
\begin{equation}
\label{Tric1}\Psi^2(a, c)\;\le \;\Psi(a-1, c-1)\Psi(a+1, c+1),
\end{equation}
autre in\'egalit\'e de type Tur\'an. On a encore $\Psi^2(a,c)\sim\Psi(a-1,c-1)\Psi(a,c+1)$ en $+\infty,$ et on peut donc se poser la question de la taille de l'\'ecart dans l'in\'egalit\'e. Le th\'eor\`eme 2 dans \cite{BI} montre \`a cet \'egard que
\begin{equation}
\label{Tric2}\Psi(a-1,c-1)\Psi(a+1,c+1)\, - \,\Psi^2(a,c)\; \le \; -\frac{1}{c}\,\Psi^2(a,c)
\end{equation}
pour tout $c<0<a,$ en utilisant \`a nouveau par la repr\'esentation int\'egrale d'Ismail et Kelker. Une cons\'equence imm\'ediate du corollaire 1 et de la transformation de Kummer est l'in\'egalit\'e suivante, dont le domaine de validit\'e est en partie plus large que celui de (\ref{Tric2}) et qui donne aussi une estimation plus pr\'ecise sur $\{x \ge -c > 0\}.$ 

\begin{Corrd} 
On a
\begin{equation}
\label{Tric3}
\Psi(a-1,c-1)\Psi(a+1,c+1)\, - \,\Psi^2(a,c)\; \le \; \frac{1}{x}\,\Psi^2(a,c)\Psi(a+1,c+1)
\end{equation}
sur $(0,+\infty),$ pour tout $a > 1$ et $c < 1+a.$ 
\end{Corrd}

Il est envisageable que d'autres in\'egalit\'es de type Tur\'an pour d'autres fonctions sp\'eciales d\'ecoulent d'arguments simples m\^elant produits de variables al\'eatoires et positivit\'e compl\`ete. Ceci pourrait faire l'objet de recherches plus pouss\'ees. \\

\noindent
{\bf Remerciements.} Je remercie G\'erard Letac pour ses commentaires. Ce travail a b\'en\'efici\'e d'une aide de l'Agence Nationale de la Recherche portant la r\'ef\'erence ANR-09-BLAN-0084-01.

\end{document}